\documentclass[11pt]{article}
\usepackage{graphicx,epsfig,amssymb,cite}
\usepackage{amsmath}

\newtheorem{theorem}{Theorem}
\newtheorem{remark}{Remark}
\newtheorem{lemma}{Lemma}

\newtheorem{corollary}{Corollary}

\newtheorem{definition}{Definition}

\date{}

\begin{document}
\title{Robust $H_\infty$ Filtering for Nonlinear Discrete-time Stochastic Systems}
\author{ Tianliang~Zhang ${}^1$,
 Feiqi~Deng ${}^1$
 \thanks{
  Email: t\_lzhang@163.com(T. Zhang), aufqdeng@scut.edu.cn(F. Deng), w\_hzhang@163.com(W. Zhang)
 }\ and Weihai~Zhang ${}^2$ \\
%EndAName
${}^1$ {\small School of Automation Science and Engineering,} \\
 {\small South China University of Technology,
Guangzhou 510640,   P. R. China}\\
${}^2$ {\small College of Information and Electrical Engineering,} \\
{\small Shandong University of Science and Technology,}  {\small
Qingdao  266510, P. R. China}\\}
\maketitle

{\bf Abstract-} {\hspace{0.13in}  This paper mainly discusses the  $H_{\infty}$
filtering of general   nonlinear  discrete time-varying  stochastic
systems. A nonlinear discrete-time stochastic bounded real lemma
(SBRL) is firstly obtained by means of the smoothness of the
conditional mathematical expectation, and then, based on the given
SBRL and a stochastic LaSalle-type theorem, a sufficient condition
for the existence of the $H_\infty$ filtering of general nonlinear
discrete time-varying stochastic systems is presented via  a new
introduced Hamilton-Jacobi inequality (HJI), which is easily
verified. When the worst-case disturbance $\{v^*_k\}_{k\in {\mathcal
N}}$ is considered,  the suboptimal $H_2/H_\infty$ filtering is
studied. Two examples including a practical engineering example show
the effectiveness of our  main results. }


\section{Introduction}

\hspace{0.13in}  $H_{\infty}$ control theory was initially
formulated by G. Zames \cite{Zames} in the early 1980's for linear
time-invariant systems, which has been one of the most important
control approaches in the presence of external disturbances. Because
in engineering practice, the system state is not always available,
how to estimate the   unavailable state variable or a linear
combination of the state variable  from  the  measurement output is
an important issue of modern control theory. When the system noise
is stationary Gaussian white noise, Kalman filtering has been shown
to be one of the most celebrated estimation methods. However, in
practical applications, we may not be able to  accurately know the
statistical properties of external disturbances. In this case, one
has to turn to Robust $H_\infty$ filter.
 Robust $H_\infty$ filter requires one to design a filter such that the $l_2$-gain from the exogenous disturbance to the
 estimated error is less than a prescribed level $\gamma>0$.
In contrast with the well-known Kalman filtering, one of the main
advantages of $H_{\infty}$ filtering is that it is not necessary to
know exactly the statistical properties of the external disturbance
but only assumes the external disturbance to have bounded energy
\cite{zhangfilter}. We refer the reader to
\cite{bschen_2008_2,Chen2017,Deng,Fault2016,WangT2017,Yan2016} for
practical applications of $H_{\infty}$ filtering in signal
processing and sensor networks.

Stochastic $H_\infty$ control of  linear continuous-time It\^o
stochastic systems seems to start from the well-known works
\cite{Hinrichsen,ug}. After then, based on the stochastic bounded
real lemma of \cite{Hinrichsen}, full- and reduced-order robust
$H_{\infty}$ estimation problems for stationary continuous-time
linear stochastic It\^o  systems were discussed in \cite{Gershon}
and \cite{Xu2002}, respectively. All the above works are limited to
the linear stationary stochastic systems. We refer the reader to the
monographs \cite{v.dragan_2005,Petersen,book} for the early
development in the $H_{\infty}$ control theory of linear It\^o
systems.
 By means of
completing squares and stochastic dynamic programming principle, the
state-feedback $H_{\infty}$ control  was  extensively investigated
in \cite{zhangsiam}  for affine stochastic It\^o systems. Based on
the stochastic bounded real lemma given in \cite{zhangsiam}, the
reference \cite{zhangfilter} solved the  the nonlinear stochastic
$H_\infty$ filter design of nonlinear affine It\^o systems by
solving a second-order Hamilton-Jacobi  inequality (HJI).

As said by J. P. LaSalle \cite{LaSalle}, ``Today there is more and
more reason for studying difference equations systematically. They
are in their own right important mathematical models.''.  Along the
development of  computer technique, it is expected that the study on
discrete-time systems will become more and more important
\cite{Shaked2006,dragan_2010,Bouhtouri1999,LinIEEE,zhao-deng-2014_1,Lin2016}.
$H_\infty$ control of linear discrete-time stochastic systems with
multiplicative noise was initiated by \cite{Bouhtouri1999}, and then
generalized to nonlinear discrete stochastic systems
\cite{Shaked2006}. In \cite{Niu2010}, $H_\infty$ filtering of
discrete fuzzy stochastic systems with sensor nonlinearities was
studied. In \cite{LiShi}, $H_{\infty}$ filtering for a class of
nonlinear discrete-time stochastic systems with uncertainties and
random Markovian delays was investigated. Although, the $H_\infty$
control and filtering of continuous-time It\^o systems have been
solved in \cite{zhangsiam} and \cite{zhangfilter}, respectively.
However, for a general nonlinear discrete  stochastic system, its
$H_{\infty}$ control and filtering problems  seem   more complicated
than continuous-time It\^o systems.  The main reason lies in that
discrete nonlinear stochastic systems do not have  an infinitesimal
generator $LV (x)$ as in It\^o systems, which is a useful tool in
completing squares \cite{Mao}. In \cite{Shaked2006}, the  nonlinear
discrete-time $H_\infty$ control was discussed based on an HJI,
where the HJI
 depends on the supremum of a conditional mathematical
expectation function, which is not easily verified. Generally
speaking, the technique of completing squares used in
\cite{zhangsiam}  becomes invalid for general discrete-time
nonlinear systems. Moreover, due to adaptiveness requirement, the
method of  Taylor's series expansion is not applicable as done in
deterministic nonlinear systems \cite{weilin-1996}; see \cite{book}.

As summarized above,   how to give practical criteria for general
discrete stochastic $H_{\infty}$ control and filtering  that do not
depend on the mathematical expectation of the state trajectory  is a
challenging work.
 In \cite{Lin2016}, by  Doob's  super-martingale theory, a LaSalle-type stability
theorem   was  established.  A new method   based on convex analysis
was introduced to solve the $H_{\infty}$ control  of general
discrete-time nonlinear stochastic systems in \cite{LinIEEE}. In
\cite{LinIEEE}, the  Lyapunov function is  selected as a convex
function, which help separate the state $x_k$ from the coupling of
$x_k$ and the unknown exogenous disturbance $v_k$.

In this paper, our main goal is to deal with $H_{\infty}$ filtering
 for general  nonlinear  discrete stochastic systems. Firstly, by applying the smoothness of the conditional mathematical
expectation, for general  discrete time-varying  nonlinear
time-varying stochastic systems, a stochastic bounded real lemma
(SBRL) on external stability  is given based on a new introduced
HJI,  where the HJI  does not depend on the mathematical
expectations of the state and external disturbance.  Secondly, a
sufficient condition for  the existence  of   $H_{\infty}$ filtering
for general nonlinear discrete time-varying   stochastic systems has
been presented based on our newly developed stochastic LaSalle's
invariant principle \cite{Lin2016} and SBRL.  As corollaries,
$H_\infty$ filtering problems   of
 nonlinear stochastic time-invariant systems and  affine nonlinear stochastic systems are  discussed.
 Thirdly,  we also discuss  the suboptimal $H_2/H_{\infty}$ filtering problems of
 nonlinear stochastic systems and  linear stochastic systems under worst-case disturbance $\{v^*_k\}_{k\in {\mathcal N}}$.
 In particular, for linear stochastic systems,
 we prove that a  desired  suboptimal  $H_2/H_{\infty}$ filtering can be constructed by solving a convex optimization problem.

 This paper is organized as follows: In section 2, some preliminaries are made, where a useful
 lemma-Lemma 2.3
 on stability in probability is obtained. Section 3 is concerned
 about the general nonlinear $H_\infty$ filtering, and section 4 is
 about the suboptimal $H_2/H_\infty$ filtering. In section 5, we
 present two examples, one is a numerical example, but the other one
 is a practical vehicle roll example, to illustrate the validity of
 our main results.

For convenience, the notations adopted in this paper are as follows:

$M'$: the transpose of the matrix $M$ or vector $M$; $M>0$ ($M<0$):
the matrix $M$ is a positive definite (negative definite) real
symmetric matrix; $I_n$: $n\times n$ identity matrix; ${\mathcal
R}^n$: the $n$-dimensional real Euclidean vector space with the
standard $2$-norm $\|x\|$; ${\mathcal R}^{n\times m}$: the space of
$n\times m$  real  matrices. ${\mathcal N}:=\{0,1,2,\cdots,\}$;
 $l^2_{\infty}(\Omega, {\mathcal F},
{\{{\mathcal F}_k\}}_{k\in \mathcal N}, \mathcal{R}^{n_v})$: the
space of stochastic ${\mathcal F}_k$-adapted   sequence
$\{v_k\}_{k\in \mathcal N}$ with the norm
$$
\|v\|_{l^2_{\infty}}=\sqrt{E\left[\sum_{k=0}^{\infty}\|v_k\|^2\right]}<\infty.
$$
A function $f(x)$ is called a positive function, if $f(x)>0$ for any
$x\neq0$, and $f(0)=0$;  $\mathcal {K}$: the family of all strictly
increasing continuous  positive  functions $\mu(\cdot)$.

\section{Preliminaries}
Consider the following general discrete time-varying nonlinear
stochastic system
\begin{equation}\label{2017filter-1}
\begin{cases}
 \  x_{k+1}=f_k(x_k,w_k,v_k), \ f_k(0,0,0)\equiv0,\\
 \  y_k=g_k(x_k,v_k), \ g_k(0,0)\equiv0,\\
    z_k=m_k(x_k,v_k), \  m_k(0,0)\equiv0,
   \end{cases}
\end{equation}
where $x_{k}\in \mathcal{R}^{n}$ is the $n$-dimensional state
vector, $y_k\in \mathcal{R}^{n_y}$ is the $n_y$-dimensional
measurement output, $v_k\in \mathcal{R}^{n_v}$ stands for the
exogenous disturbance signal with $\{v_k\}_{k\in \mathcal N} \in
l^2_{\infty}(\Omega, {\mathcal F}, {\{{\mathcal F}_k\}}_{k\in
\mathcal N}, \mathcal{R}^{n_v})$, $z_k\in \mathcal{R}^{n_z}$ is
called the regulated output, which  is the combination of the  state
and the exogenous disturbance to be estimated,  and
$\{w_k\}_{k\in\mathcal N}$ is a sequence of independent
$n_w$-dimensional random variables defined on the complete  filtered
probability space $(\Omega, {\mathcal F}, {\{{\mathcal F}_k\}}_{k\in
\mathcal N}, {\mathcal P})$, where ${\mathcal
F}_k=\sigma(w_0,w_1,\cdots,w_{k-1})$, ${\mathcal
F}_0=\{\phi,\Omega\}$.  $f_k: \mathcal{R}^{n}   \times
\mathcal{R}^{n_w} \times \mathcal{R}^{n_v}\mapsto \mathcal{R}^{n}$ ,
$g_k:\mathcal{R}^{n} \times\mathcal{R}^{n_v}
\mapsto\mathcal{R}^{n_y}$ and $m_k:\mathcal{R}^{n}
\times\mathcal{R}^{n_v}\mapsto \mathcal{R}^{n_z} $ are continuous
vector-valued functions.  $\{x_k^{s,x,v}\}_{k\in\mathcal N}$ denotes
the solution sequence of system (\ref{2017filter-1}) with the
initial state $x\in {\mathcal R}^{n}$ starting at the initial time
$k=s$ under the exogenous disturbance $\{v_k\}_{k\in \mathcal N}$.
Similarly, $\{y_k^{s,x,v}\}_{k\in\mathcal N}$ and
$\{z_k^{s,x,v}\}_{k\in\mathcal N}$ can also be defined.
%%%%%%%%%%%%%%%%%%%%%%%

In what follows, we construct the following  filter for the
estimation of $z_k$:
\begin{equation}\label{2017filter-2}
\begin{cases}
 \  \hat{x}_{k+1}=\hat{f}_k(\hat{x}_k)+\hat{g}_k(y_k), \  \hat{f}_k(0)=\hat{g}_k(0)=0,\\
    \hat{z}_k=\hat{m}_k(\hat{x}_k),  \ \hat{m}_k(0)=0, \ \hat{x}_0=0,
   \end{cases}
\end{equation}
where $\hat{x}_{k} \in \mathcal {R}^{n_{\hat{x}}}$ is the estimated
value of $x_k$,  $\hat{f}_k$, $\hat{g}_k$, and $\hat{m}_k$ are
filter parameters  to be determined. $\hat{z}_k \in \mathcal
{R}^{n_{z}}$ is the estimated value of $z_k$ in system
(\ref{2017filter-1}). Let
$\eta_k=\left[\begin{array}{cc}x_k\\\hat{x}_k\end{array}\right]$,
and $\tilde{z}_k\triangleq z_k-\hat{z}_k$, where $\tilde{z}_k$
denotes the estimation error of the regulated output, then we get
the following augmented system:
\begin{equation}\label{2017filter-3}
\begin{cases}
 \  \eta_{k+1}=\tilde{f}_k(\eta_k, w_k,v_k)=\left[\begin{array}{cc}f_k(x_k,w_k,v_k)\\ \hat{f}_k(\hat{x}_k)+\hat{g}_k(y_k) \end{array}\right], \\
 y_k=g_k(x_k,v_k),\\
    \tilde{z}_k= \tilde{m}_k(\eta_k,v_k) = m_k(x_k,v_k)-\hat{m}_k(\hat{x}_k),
   \end{cases}
\end{equation}
where $\eta_{k}\in \mathcal {R}^{n_{\eta}}$.

\begin{remark}\label{remark1}
In practical engineering, only $\dim (\hat{x}_{k})\le \dim (x_k)$ is
valuable. When $\dim (\hat{x}_{k})=\dim (x_k)$, (\ref{2017filter-2})
is called a full-order filter; when $\dim (\hat{x}_{k})<\dim (x_k)$,
(\ref{2017filter-2}) is called a reduced-order filter; see
\cite{Xu2002}.
\end{remark}

%%%%%%%%%%%%%%%%%%%%%%%%%%%%%%%%%%%%%%%%%%%%%%%%%%%%%%%%%%%%%%%%%%%%%%%%%%%%%%%%%%%%%%%%%%%%
\begin{definition}{(\bf Internal stability)}\label{2017filter-def1}
System (\ref{2017filter-3}) is said to be internally stable, if when
$v=\{v_k\}_{k\in {\mathcal N}}\equiv0$, the zero solution of
\begin{equation}\label{eqvgvh}
\eta_{k+1}=\tilde{f}_k(\eta_k,w_k,0)
\end{equation}
is globally  asymptotically stable in probability. In other word,
for any $\varepsilon>0$, there have
\begin{equation}\label{2017filter-rrt1}
\lim_{\eta_0\rightarrow0}\mathcal {P}\left(\sup_{k\geq
0}\|\eta_k^{0,\eta_0,0}\|>\varepsilon\right)=0
\end{equation}
 and
\begin{equation}\label{2017filter-rrt2}
{\mathcal P} \{\lim_{k\to \infty} \eta_k^{0,\eta_0,0}=0\}=1.
\end{equation}
When only (\ref{2017filter-rrt1}) holds, system (\ref{eqvgvh}) is
said to be stable in probability.
\end{definition}
%%%%%%%%%%%%%%%%%%%%%%%%%%%%%%%%%%%%%%%%%%%
\begin{definition}\textbf {(External stability)}\label{2017filter-def2}
For any given positive real number $\gamma>0$, system
(\ref{2017filter-3}) is called externally stable, if for  any
nonzero $\{v_k\}_{k\in \mathcal N} \in l^2_{\infty}(\Omega,
{\mathcal F}, {\{{\mathcal F}_k\}}_{k\in \mathcal N},
\mathcal{R}^{n_v})$  and zero initiate state $\eta_0=0$, we have
\begin{equation}\label{2017filter-eqw}
\|\tilde{z}\|^2_{l^2_{\infty}} \leq \gamma^2\|v\|^2_{l^2_{\infty}}.
\end{equation}
\end{definition}
%%%%%%%%%%%%%%%%%%%%%%%%%%%%%%%%%%%%%%%%%%%%%%%%%%%%%%%

\begin{remark}
Define an operator called the perturbation operator of system
(\ref{2017filter-3}) as follows:
\begin{eqnarray*}
\ \ \ \ \ \ \ \ \ \ \ \ \ \ \mathcal {L}:&&v \in
l^2_{\infty}(\Omega, {\mathcal F}, {\{{\mathcal F}_k\}}_{k\in
\mathcal N}, \mathcal{R}^{n_v})  \\ &&\mapsto \tilde{z}\in
l^2_{\infty}(\Omega, {\mathcal F}, {\{{\mathcal F}_k\}}_{k\in
\mathcal N}, \mathcal{R}^{n_z})
\end{eqnarray*}
with the norm of ${\mathcal L}$  defined by
\begin{equation}\label{2017filter-wq}
\|\mathcal {L}\|=\sup_{\eta_0=0, 0\neq v\in l^2_{\infty}(\Omega,
{\mathcal F}, {\{{\mathcal F}_k\}}_{k\in \mathcal N},
\mathcal{R}^{n_v})}\frac{\|\tilde{z}\|_{l^2_{\infty}}}{
\|v\|_{l^2_{\infty}}},
\end{equation}
then the inequality (\ref{2017filter-eqw}) can be rewritten as
$\|\mathcal {L}\|\leq\gamma$.  $\|\mathcal {L}\|$ means the worst
case effect from the stochastic disturbance $\{v_k\}_{k\in\mathcal
{N}}$ to the controlled output $\{\tilde{z}_k\}_{k\in\mathcal {N}}$.
Therefore, in order to determine whether the system
(\ref{2017filter-3}) is externally  stable, it is important to find
a way to determine or estimate the norm $\|\mathcal {L}\|$.
\end{remark}

The nonlinear stochastic $H_{\infty}$ filtering can be stated as
follows:

\begin{definition}\textbf {(Nonlinear stochastic $H_{\infty}$ filtering)}
\label{2017filter-4} Find the filter parameters
$\hat{f}_k(\hat{x})$, $\hat{g}_k(y)$, and $\hat{m}_k(\hat{x}_k)$
such that
\begin{description}
  \item[(i)] The augmented system (\ref{2017filter-3}) is internally stable,  i.e., when $v_k=0$,
  ${k\in\mathcal{N}}$, system (\ref{2017filter-3}) is globally  asymptotically stable in probability.
  \item[(ii)] The augmented system (\ref{2017filter-3}) is externally stable,  i.e.,  the norm of the  perturbation operator
      \begin{equation*}
\|\mathcal {L}\|\leq\gamma,
\end{equation*}
where $\gamma>0$ is the given disturbance attenuation level.
\end{description}
\end{definition}

\begin{definition} \label{eqbhbn} \cite{Has'minskii,Mao} We consider a continuous function $V_k(x):=V(k,x)$ defined on ${\mathcal {N}} \times \mathcal{D}_r$
with $\mathcal{D}_r:=\{x:\|x\|\le r\}$. Let ${\mathcal K}$ be the
family of all continuous strictly increasing  functions $\mu(\cdot):
{\mathcal R}^+\mapsto {\mathcal R}^+$, such that $\mu(0)=0$ and
$\mu(t)>0$ for any $t>0$.

\begin{itemize}

\item [(i)] $\{V_k(x)\}_{k\in {\mathcal {N}}}$ is a positive definite function sequence on ${\mathcal {N}} \times \mathcal{D}_r$ in
the sense of Lyapunov if $V_k(0)\equiv 0$ for $k\in {\mathcal {N}}$,
and there exists $\mu\in {\mathcal K}$, such that
    $$
    V_k(x)\ge \mu(\|x\|), \ \ \forall (k,x)\in {\mathcal {N}} \times \mathcal{D}_r.
    $$

\item [(ii)] $\{V_k(x)\}_{k\in {\mathcal {N}}}$ is said to be  radially
unbounded if
\begin{equation}\label{radial unbounded condition}
\liminf_{\|x\|\to \infty} \inf_{k\in {\mathcal {N}}} V_k(x)=\infty.
%\inf_{k>0}V(x(k))\rightarrow\infty, \ \  as \ \ ||x||\rightarrow\infty.
\end{equation}
\end{itemize}
\end{definition}

We call $v^*:=\{v_k=v_k^*\}_{k\in {\mathcal N}}$ as the worst-case
disturbance sequence, if
\begin{equation*}
v^*=\arg min
\{J_{\infty}(v):=E\sum_{k=0}^{\infty}(\gamma^2\|v_k\|^2-\|\tilde{z}_k\|^2)\},
\end{equation*}
see \cite{book}.

The following property of the  conditional mathematical expectation
plays an important role in this paper.

\begin{lemma}[Theorem 6.4 of \cite{Kallenberg2002}]\label{2017filter-5}
If ${\mathcal R}^n$-valued random variable $\alpha$ is independent
of the $\sigma$-field $\mathcal {G}$, and ${\mathcal R}^d$-valued
random variable $\beta$ is $\mathcal {G}$-measurable, then for any
bounded or nonnegative function $f:{\mathcal R}^n\times {\mathcal
R}^d\mapsto {\mathcal R}$,
\begin{equation*}
E[f(\alpha,\beta)|\mathcal {G} ]=E[f(\alpha,x)]_{x=\beta}\ \ \ a.s.
\end{equation*}
holds.
\end{lemma}
For a Lyapunov  function sequence $\{V_k(\eta)\}_{k\in{\mathcal
N}}$, $V_k: {\mathcal R}^{n_{\eta}} \mapsto {\mathcal R}^+$, because
$w_k$ is independent of $\mathcal {F}_{k}$,  $x_k$,  ${\hat x}_k$
and $v_k$ are $\mathcal {F}_{k}$-measurable, from Lemma
\ref{2017filter-5}, we have
\begin{eqnarray*}
E[V_{k+1}(\eta_{k+1})|\mathcal {F}_{k}]&=&E[V_{k+1}(\tilde{f}_k(\eta_k, w_k,v_k))|\mathcal {F}_{k}]\nonumber\\
&=&E[V_{k+1}(\tilde{f}_k(\eta, w_k,v))]_{\eta=\eta_k, v=v_k}.
\end{eqnarray*}
We define
\begin{eqnarray*}
&&\Delta_{v} V_k(\eta)\\&:=&E[V_{k+1}(\tilde{f}_k(\eta,
w_k,v))]-V_k(\eta),\ \ \ \eta\in {\mathcal R}^{n_{\eta}}, v \in
{\mathcal R}^{n_{v}}
\end{eqnarray*}
as the difference operator of the function sequence
$\{V_k\}_{k\in{\mathcal N}}$, and set $H_k(\eta,v):=\Delta_{v}
V_k(\eta)+\|\tilde{m}_k(\eta,v)\|^2$.

The following lemma is referred to as a SBRL, which gives a
sufficient condition for the external stability.
\begin{lemma}\label{2017filter-lemma2}
For a given $\gamma>0$, if there exists a positive definite Lyapunov
function sequence $\{V_k(\eta)\}_{k\in{\mathcal N}}$ satisfies the
following HJI inequality
\begin{equation}\label{2017filter-rree}
H_k(\eta,v)-\gamma^2\|v\|^2\leq0, \   k\in {\mathcal N},
\end{equation}
then  system (\ref{2017filter-3}) is  externally stable.
\end{lemma}
\textbf {Proof}: Let $\{\eta_k\}_{k\in {\mathcal N}}$ be the
solution of system (\ref{2017filter-3}), we have
\begin{align}\label{rebjh}
&E[V_{k+1}(\eta_{k+1})]-E[V_k(\eta_k)]\nonumber\\
=&E\{E[V_{k+1}(\eta_{k+1})|{\mathcal F}_k]\}-E[V_k(\eta_k)]\nonumber\\
=&E\{E[V_{k+1}(\tilde{f}_k(\eta_k, w_k,v_k)|{\mathcal
F}_k]-V_k(\eta_k)\}.
\end{align}
By Lemma~\ref{2017filter-5}, we have
\begin{align*}
&E\{E[V_{k+1}(\tilde{f}_k(\eta_k, w_k,v_k)|{\mathcal
F}_k]-V_k(\eta_k)\}\\
=&E\{E[V_{k+1}(\tilde{f}_k(\eta,
w_k,v))]_{\eta=\eta_k,v=v_k}-V_k(\eta_k)\}\\
=&E\{[\Delta_v V_k(\eta)]_{\eta=\eta_k, v=v_k}\},
\end{align*}
which together with condition (\ref{2017filter-rree}), it follows
from (\ref{rebjh}) that
\begin{align}\label{2017filter-rsgsgfe}
&E[V_{k+1}(\eta_{k+1})]-E[V_k(\eta_k)]\nonumber\\ \leq& \gamma^2
E\|v_k\|^2-E\|\tilde{m}_k(\eta_k,v_k)\|^2.
\end{align}
Taking a summation on both sides of the above inequality from $k=0$
to $k=N$, we can get that
\begin{align}\label{2017filter-rsgs}
&E[V_{N+1}(\eta_{N+1})]-E[V_0(\eta_0)]\nonumber\\
\leq& \gamma^2\sum_{k=0}^{N}
E\|v_k\|^2-\sum_{k=0}^{N}E\|\tilde{m}_k(\eta,v_k)\|^2.
\end{align}
Note that $V_k(\eta)>0$ for $\eta\neq0$ due to the positivity of the
function sequence $\{V_k\}_{k\in{\mathcal N}}$ and $\eta_0=0$, then
(\ref{2017filter-rsgs}) leads to
\begin{align*}
\sum_{k=0}^{N}E\|\tilde{m}_k(\eta_k,v_k)\|^2=&
\sum_{k=0}^{N}E\|\tilde{z}_k\|^2\\ \leq& \gamma^2
\sum_{k=0}^{N}E\|v_k\|^2.
\end{align*}
Let $N\rightarrow\infty$, we obtain
$$
\sum_{k=0}^{\infty}E\|\tilde{z}_k\|^2\leq \gamma^2
\sum_{k=0}^{\infty}E\|v_k\|^2
$$
for  any $\{v_k\}_{k\in \mathcal N} \in l^2_{\infty}(\Omega,
{\mathcal F}, {\{{\mathcal F}_k\}}_{k\in \mathcal N},
\mathcal{R}^{n_v})$. The proof is completed according to
Definition~\ref{2017filter-def2}. $\square$

\begin{lemma}\label{2017filter-lemmastablefffaa} Suppose there exist a  positive definite Lyapunov
 function sequence $\{V_k(\eta)\}_{k\in{\mathcal N}}$  such that
\begin{equation}
\Delta_{v=0} V_k(\eta)\leq 0, \  \forall (\eta,k) \in {\mathcal D}_h
\times {\mathcal N}, \label{dd2017filter-L1}
\end{equation}
 then the  stochastic difference
equation
$$
\eta_{k+1}=\tilde{f}_k(\eta_k, w_k,0),  \ \eta_{0} \in {\mathcal
R}^\eta
$$
is stable in probability.
\end{lemma}

\textbf {Proof}:  By Definition~\ref{eqbhbn}-(i), there exists
$\mu(\cdot)\in {\mathcal K}$ such that
$$
V_k(\eta)\ge \mu(\|\eta\|), \forall (\eta,k)\in {\mathcal D}_r
\times {\mathcal N}.
$$
For any $\epsilon\in (0,1)$ and $0<r<h$, because $V_k(0)=0$, $k\in
{\mathcal N}$, $V_k(\eta)$ is continuous with respect to $\eta$, we
can find a positive constant $\delta>0$, $\delta<r$ such that
\begin{equation}\label{eq vhgv}
\frac 1\epsilon \sup_{\eta\in {\mathcal D}_\delta} V_0(\eta)\le
\mu(r).
\end{equation}
For any fixed $\eta_0\in {\mathcal D}_\delta$, define
$$
\tau_r=\max\{k:\|\eta_k\|<r\}.
$$
Under the condition of $\Delta_{v=0} V_k(\eta)\leq 0$, by
Lemma~\ref{2017filter-5},  for any $\eta_i\in  {\mathcal D}_r$, we
have
\begin{align*}
 &E[V_{i+1}(\eta_{i+1})|{\mathcal F}_i]-V_i(\eta_{i})\\=&E[V_{i+1}(\tilde{f}_i(\eta_i,\omega_i,0))-V_i(\eta_{i})|{\mathcal F}_i]\\
 =&E[V_{i+1}(\tilde{f}_i(\eta,\omega_i,0))-V_i(\eta)]\big|_{\eta=\eta_i}\\
 =& \Delta_{v=0} V_i(\eta)\big|_{\eta=\eta_i}\le 0, \  a.s..
\end{align*}
Taking a summation from $i=0$ to $\tau_r \wedge k$, it follows that
\begin{equation*}
\sum_{i=0}^{\tau_r \wedge k} [E[V_{i+1}(\eta_{i+1})|{\mathcal
F}_i]-V_i(\eta_{i})]\le 0.
\end{equation*}
Taking a mathematical expectation operator in above, we have
\begin{equation}\label{eqbhbdc}
E[V_{\tau_r \wedge k+1}(\eta_{\tau_r \wedge k+1})]\le V_0(\eta_0).
\end{equation}
Because
\begin{align}
&E[V_{\tau_r \wedge k+1}(\eta_{\tau_r \wedge
k+1})]\nonumber\\
=&E[I_{\{{\tau_r}<k\}}
V_{\tau_r+1}(\eta_{\tau_r+1})]+E[I_{\{{\tau_r}\ge k\}} V_{k+1}(\eta_{k+1})]\nonumber\\
\ge& E[I_{\{{\tau_r}<k\}} V_{\tau_r+1}(\eta_{\tau_r+1})]\nonumber\\
\ge & {\mathcal P}(\tau_r<k)\mu(r).\label{eqdeda}
\end{align}
Combing (\ref{eq vhgv}), (\ref{eqbhbdc}) and (\ref{eqdeda}) leads to
$$
{\mathcal P}(\tau_r<k)\le \frac {V_0(\eta_0)}{\mu(r)}\le \epsilon.
$$
Let $k\to\infty$, then  ${\mathcal P}(\tau_r<\infty)\le
\varepsilon$, or equivalently, ${\mathcal P}(\|\eta_k\|<r, \forall
k\in {\mathcal N})>1-\varepsilon$. This theorem is proved. $\square$

The following is the so-called LaSalle-type theorem, which cites
from  Theorem 3.1 of  \cite{Lin2016}.

\begin{lemma}\label{2017filter-lemmastablefff} Suppose there exist a  radially unbounded positive  Lyapunov  function
sequence $\{V_k(\eta)\}_{k\in{\mathcal N}}$, a deterministic
real-valued sequence $\{\gamma_k\geq0\}_{k\in{\mathcal N}}$, and a
nonnegative function $W: {\mathcal R}^{n_{\eta}}\mapsto{\mathcal
R}^+$, satisfying
\begin{equation}
\Delta_{v=0} V_k(\eta)\leq \gamma_k-W(\eta), \forall \eta \in
{\mathcal R}^{n_{\eta}},  k\in {\mathcal N}, \label{2017filter-L1}
\end{equation}
\begin{equation}
\sum_{k=0}^{\infty}\gamma_k<\infty.\label{2017filter-L2}
\end{equation}
Let  $\eta_{k}$  be the solution of
$$
\eta_{k+1}=\tilde{f}_k(\eta_k, w_k,0), \ k\in {\mathcal N},
$$
then
$$
\lim_{k\rightarrow\infty} V_k(\eta_k) \ \ \mbox{exists and is finite
almost surely},
$$
and
$$
\lim_{k\rightarrow\infty}W(\eta_k)=0, a.s..
$$
\end{lemma}

\section{General nonlinear  $H_{\infty}$ filtering}

In this section, we will discuss the $H_\infty$ filtering design
problem for both nonlinear time-varying and time-invariant
stochastic systems.

\begin{theorem}\label{2017filter-th5}
For a given disturbance attenuation level $\gamma>0$, suppose that
 there exist a positive definite   radially unbounded  Lyapunov
sequence $\{V_k\}_{k\in{\mathcal N}}$, and a positive radially
unbounded function $W: {\mathcal R}^{n_{\eta}}\mapsto{\mathcal
R}^+$, such that  (\ref{2017filter-rree}) and
\begin{equation}\label{eqvhvh}
\|{\tilde m_k}(\eta,0)\|^2\ge W(\eta)
\end{equation}
hold, then system (\ref{2017filter-2}) is a desired $H_{\infty}$
filter for system (\ref{2017filter-1}).
\end{theorem}
\textbf{Proof:} We first show that  the augmented system
(\ref{2017filter-3}) is internally stable. By
(\ref{2017filter-rree}) and the definition of $H_k(\eta,v)$, we have
$$
\Delta_{v=0} V_k(\eta)+\|\tilde{m}_k(\eta,0)\|^2\le 0
$$
or equivalently,
\begin{equation}\label{eqgghg}
\Delta_{v=0} V_k(\eta)\le -\|\tilde{m}_k(\eta,0)\|^2.
\end{equation}
By (\ref{eqvhvh}), (\ref{eqgghg}) gives that
$$
\Delta_{v=0} V_k(\eta)\le -W(\eta),
$$
which, according to Lemma~\ref{2017filter-lemmastablefffaa}, implies
that system (\ref{2017filter-3}) is stable in probability when
$v=0$. In addition, by Lemma~\ref{2017filter-lemmastablefff}, we
have $\lim_{k\to\infty} W(\eta_k)=0, a.s..$  In view of $W(\eta)$
being  a positive radially unbounded  function, hence, the set of
all limit points $G^*(\eta)$ only contains a zero point, that is,
$$
G^*(\eta)=\{\eta: W(\eta)=0\}=\{0\}.
$$
So from  $\lim_{k\to\infty} W(\eta_k)=0$, a.s.,  it follows that
$$
{\mathcal P}\{\lim_{k\to\infty} \eta_k=0\}=0.
$$
The internal stability is proved. While the external stability  is
obtained by Lemma~\ref{2017filter-lemma2}. The proof of this theorem
is completed.  $\square$

In particular, for the following  discrete time-invariant nonlinear
stochastic system
\begin{equation}\label{2017filter-new1a}
\begin{cases}
 \  x_{k+1}=f(x_k,w_k,v_k),\ f(0,0,0)=0, \\
 \  y_k=g(x_k,v_k),\ g(0,0)=0,\\
    z_k=m(x_k,v_k),\ m(0,0)=0,
   \end{cases}
\end{equation}
the filter equation is often taken as
\begin{equation}\label{2017filter-new2}
\begin{cases}
 \  \hat{x}_{k+1}=\hat{f}(\hat{x}_k)+\hat{g}(y_k),\  \hat{f}(0)=\hat{g}(0)=0, \\
    \hat{z}_k=\hat{m}(\hat{x}_k),\ \hat{m}(0)=0,\ \hat{x}_0=0.
   \end{cases}
\end{equation}
In this case, the augmented system can be rewritten as
\begin{equation}\label{2017filter-new3}
\begin{cases}
 \  \eta_{k+1}=\tilde{f}(\eta_k, w_k,v_k)=\left[\begin{array}{cc}f(x_k,w_k,v_k)\\ \hat{f}(\hat{x}_k)+\hat{g}(y_k) \end{array}\right], \\
    \tilde{z}_k= \tilde{m}(\eta_k,v_k) = m(x_k,v_k)-\hat{m}(\hat{x}_k).
   \end{cases}
\end{equation}
We choose a common positive  Lyapunov function $V: {\mathcal
R}^{n_{\eta}}\mapsto{\mathcal R}^+$, and for any $\eta\in {\mathcal
R}^{n_{\eta}}, v\in {\mathcal R}^{n_{v}}$, we write
\begin{eqnarray}\label{2017filter-caxa}
\Delta_{v,k} V(\eta):=E[V(\tilde{f}(\eta, w_k,v))]-V(\eta),  \
\forall k\in\mathcal {N}, \label{2017filter-caxa1}
\end{eqnarray}
 and
reset $H_k(\eta,v):=\Delta_{v,k} V(\eta)+\|\tilde{m}(\eta,v)\|^2$.
Theorem~\ref{2017filter-th5} immediately yields the following
corollary:

\begin{corollary}\label{2017filter-th1}
For a given  disturbance attenuation level $\gamma>0$. Suppose there
exist a positive  radially unbounded  Lyapunov function $V(\eta)$,
and a positive radially unbounded  function $W: {\mathcal
R}^{n_{\eta}}\mapsto{\mathcal R}^+$, such that
(\ref{2017filter-rree}) and
\begin{align}
\|\tilde{m}(\eta,0)\|^2\geq  W(\eta)\label{2017filter-eew2a}
\end{align}
hold, then system (\ref{2017filter-new2}) is a desired $H_\infty$
filter for system (\ref{2017filter-new1a}).
\end{corollary}

\begin{remark} If $\{w_k\}_{k\in \mathcal {N}}$ is an independently
identically distributed random variable sequence, then
\begin{equation*}
E[V(\tilde{f}(\eta, w_k,v))]=E[V(\tilde{f}(\eta, w_0,v))], \ \forall
k\in{\mathcal N},
\end{equation*}
while  (\ref{2017filter-caxa1}) and  $H_k(\eta,v)$ can be replaced
by
$$
\Delta_{v} V(\eta):=E[V(\tilde{f}(\eta, w_0,v))]-V(\eta)
$$
and
$$
H(\eta,v):=\Delta_{v} V(\eta)+\|\tilde{m}(\eta,v)\|^2,
$$
respectively.
\end{remark}

\begin{corollary}\label{2017filter-3.2}
For a  given  disturbance attenuation level $\gamma>0$, suppose
$\{w_k\}_{k\in \mathcal {N}}$ is an independently identically
distributed random variable sequence. If  there exist a positive
 radially unbounded  Lyapunov function $V(\eta)$, and a
positive radially unbounded  function $W: {\mathcal
R}^{n_{\eta}}\mapsto{\mathcal R}^+$, satisfying
(\ref{2017filter-eew2a}) and
\begin{align}
H(\eta,v)-\gamma^2\|v\|^2\leq 0,\label{2017filter-eew3a}
\end{align}
then system (\ref{2017filter-new2}) is a desired $H_\infty$ filter
for system (\ref{2017filter-new1a}).
\end{corollary}

A special case of system (\ref{2017filter-new1a}) is the following
affine nonlinear stochastic system
\begin{equation}\label{2017filter-6}
\begin{cases}
 \  x_{k+1}=f_{1}(x_k)+h_{1}(x_k)v_k+[f_{2}(x_k)+h_{2}(x_k)v_k]w_k,\\
 y_{k}=g_{1}(x_k)+g_{2}(x_k)v_k,\\
 z_k=m(x_k), \ k\in {\mathcal N},
   \end{cases}
\end{equation}
where $\{w_k\}_{k\ge 0}$ is
 a sequence of  one-dimensional independent white
noise processes. Assume that $E[w_i]=0$, $E[w_iw_j]=\delta_{ij}$,
where $\delta_{ij}$ is a Kronecker function defined by
$\delta_{ij}=0$ for $i\ne j$ while $\delta_{ij}=1$ for $i=j$.  We
take the  $H_{\infty}$ filter equation for system
(\ref{2017filter-6}) as
\begin{equation}\label{2017filter-7}
\begin{cases}
 \  \hat{x}_{k+1}=\hat{f}(\hat{x}_k)+\hat{g}(\hat{x}_k)y_k, \\
    \hat{z}_k=\hat{m}(\hat{x}_k),  \  k\in {\mathcal N}.
   \end{cases}
\end{equation}
Thus the augmented system can  be rewritten as
\begin{equation}\label{2017filter-8}
\begin{cases}
 \  \eta_{k+1}=\tilde{f}_{1}(\eta_k)+\tilde{h}_{1}(\eta_k)v_k
 +[\tilde{f}_{2}(\eta_k)+\tilde{h}_{2}(\eta_k)v_k]w_k,\\
    \tilde{z}_k= \tilde{m}(\eta_k),
   \end{cases}
\end{equation}
where
\begin{align}
&\tilde{f}_{1}(\eta_k)= \left[\begin{array}{cc}f_{1}(x_k)\\
\hat{f}(\hat{x}_k)+\hat{g}(\hat{x}_k)g_{1}(x_k)\end{array}\right], \
\tilde{f}_{2}(\eta_k)=
\left[\begin{array}{cc}f_{2}(x_k)\\0\end{array}\right], \ \nonumber\\
&\tilde{h}_{1}(\eta_k)= \left[\begin{array}{cc}h_{1}(x_k)\\
\hat{g}(\hat{x}_k)g_{2}(x_k)\end{array}\right],
\tilde{h}_{2}(\eta_k)=
\left[\begin{array}{cc}h_{2}(x_k)\\0\end{array}\right],\nonumber\\
&\tilde{m}(\eta_k) = m(x_k)-\hat{m}(\hat{x}_k).\nonumber
\end{align}

\begin{theorem}
Given the disturbance attenuation level $\gamma>0$. Suppose
$\|\tilde{m}(\eta)\|$ is a continuous positive radially unbounded
function, i.e.,
$$
\lim_{\|\eta\|\to\infty}\inf \|\tilde{m}(\eta)\|=\infty.
$$
If there exists a solution $(Q>0,{\hat f},{\hat g}, {\hat m})$
satisfying
\begin{equation}\label{2017filter-fds}
2\Theta_1+\bar{\Theta}_2-\eta'Q\eta+\|\tilde{m}(\eta)\|^2\leq0,
\end{equation}
where
\begin{equation}\label{2017filter-41}
\Theta_1=\tilde{f}'_{1}(\eta)Q\tilde{f}_{1}(\eta)
+\tilde{f}'_{2,}(\eta)Q\tilde{f}_{2}(\eta),
\end{equation}
\begin{equation}
\bar{\Theta}_2=v'(2\tilde{h}'_{1}(\eta)Q\tilde{h}_{1}(\eta)
+2\tilde{h}'_{2}(\eta)Q\tilde{h}_{2}(\eta)-\gamma^2
I)v.\label{eqvbhbv}
\end{equation}
Then, system (\ref{2017filter-7}) is the desired  $H_{\infty}$
filter of system (\ref{2017filter-6}).
\end{theorem}

\textbf{Proof}: We take $V(\eta)=\eta' Q\eta$  for all $\eta\in
{\mathcal {R}}^{n_{\eta}}$. Since $Q>0$, there holds
$$
\liminf_{\|\eta\|\to \infty}V(\eta)=\infty.
$$
For any $\eta\in {\mathcal R}^{n_{\eta}}, v \in {\mathcal
R}^{n_{v}},$
\begin{align}
&H_k(\eta,v)\nonumber\\
=&\Delta_{v,k} V(\eta)+\|\tilde{m}(\eta)\|^2\nonumber\\
=&E[V(\tilde{f}_{1}(\eta)+\tilde{h}_{1}(\eta)v
 +[\tilde{f}_{2}(\eta)+\tilde{h}_{2}(\eta)v]w_k)]-V(\eta)+\|\tilde{m}(\eta)\|^2\nonumber\\
=&\Theta_1+\Theta_2+\Theta_3-\eta'Q\eta+\|\tilde{m}(\eta)\|^2,
 \end{align}
where $\Theta_1$ is defined in (\ref{2017filter-41}),
 \begin{align*}
\Theta_2=&v'(\tilde{h}'_{1}(\eta)Q\tilde{h}_{1}(\eta)+\tilde{h}'_{2}(\eta)Q\tilde{h}_{2}(\eta))v,\\
\Theta_3=&\tilde{f}'_{1}(\eta)Q\tilde{h}_{1}(\eta)v+v'\tilde{h}'_{1}(\eta)Q\tilde{f}_{1}(\eta)
+\tilde{f}'_{2}(\eta)Q\tilde{h}_{2}(\eta)v\\
&+v'\tilde{h}'_{2}(\eta)Q\tilde{f}_{2}(\eta).
\end{align*}
Because
$$
\tilde{f}'_{1}(\eta)Q\tilde{h}_{1}(\eta)v+v'\tilde{h}'_{1}(\eta)Q\tilde{f}_{1}(\eta)\le
v'\tilde{h}'_{1}(\eta)Q\tilde{h}_{1}v+\tilde{f}'_{1}(\eta)Q\tilde{f}_{1},
$$
$$
\tilde{f}'_{2}(\eta)Q\tilde{h}_{2}(\eta)v+v'\tilde{h}'_{2}(\eta)Q\tilde{f}_{2}(\eta)\le
\tilde{f}'_{2}(\eta)Q\tilde{f}_{2}+v'\tilde{h}'_{2}(\eta)Q\tilde{h}_{2}v.
$$
Hence, according to  (\ref{2017filter-fds}),
\begin{align*}
H_k(\eta,v)-\gamma^2\|v\|^2\leq
2\Theta_1+\bar{\Theta}_2-\eta'Q\eta+\|\tilde{m}(\eta)\|^2\le 0,
\end{align*}
where $ \bar{\Theta}_2$ is defined in  (\ref{eqvbhbv}). So system
(\ref{2017filter-8}) is externally  stable by
Lemma~\ref{2017filter-lemma2}. In addition, from
$H_k(\eta,v)-\gamma^2\|v\|^2\le 0$, it yields that  $\Delta_{v=0,k}
V(\eta)\le -\|\tilde{m}(\eta)\|^2$. By
Lemma~\ref{2017filter-lemmastablefff},
$$
\lim_{k\to\infty} \|\tilde{m}(\eta_k)\|=0, a.s..
$$
Since $\|\tilde{m}(\eta_k)\|$ is a positive and radially unbounded
function, we must have $\lim_{k\to\infty} \eta_k=0$, a.s.. Hence,
system (\ref{2017filter-8}) is internally  stable.  The proof is
completed. $\square$

If  $Q$ can be partitioned as  a  block  diagonal matrix
\begin{align*}
Q=\left[
\begin{array}{cc}
Q_{1}& 0\\0 &Q_{2}
\end{array}
\right],
\end{align*}
then, the  inequality (\ref{2017filter-fds}) can be rewritten as
\begin{eqnarray*}
&&2\Theta_1+\bar{\Theta}_2-\eta'Q\eta+\|\tilde{m}(\eta)\|^2\nonumber\\
&\leq&2[f_{1}'(x)Q_{1}f_{1}(x)+(\hat{f}(\hat{x})+\hat{g}(\hat{x})g_{1}(x))'Q_{2}\nonumber\\
&&(\hat{f}(\hat{x})+\hat{g}(\hat{x})g_{1}(x))]+v'[2h'_1(x)Q_{1}h_1(x)v\nonumber\\
&&+2(\hat{g}(\hat{x})g_{2}(x))'Q_{2}
(\hat{g}(\hat{x})g_{2}(x))-\gamma^2 I]v\nonumber\\
&&-x'Q_{1}x+\hat{x}'Q_{2}\hat{x}+2\|m(x)\|^2+2\|\hat{m}(\hat{x})\|^2.\nonumber
\end{eqnarray*}
So the following corollary is obtained.

\begin{corollary}\label{2017filter-cor1}
 Consider the disturbance attenuation level $\gamma>0$.
 If there exists  the solution ($Q_1>0$, $Q_2>0$, $\hat{f}$, $\hat{g}$, $\hat{m}$) solving
\begin{eqnarray}\label{eq44}
2[h'_1(x)Q_{1}h_1(x)+(\hat{g}(\hat{x})g_{2}(x))'Q_{2}
(\hat{g}(\hat{x})g_{2}(x))]-\gamma^2 I\leq0,
\end{eqnarray}
and
\begin{eqnarray}\label{eq45}
2[f_{1}'(x)Q_{1}f_{1}(x)+(\hat{f}(\hat{x})+\hat{g}(\hat{x})g_{1}(x))'Q_{2}
(\hat{f}(\hat{x})\nonumber+\hat{g}(\hat{x})g_{1}(x))]\nonumber\\-x'Q_{1}x+\hat{x}'Q_{2}\hat{x}+2\|m(x)\|^2+2\|\hat{m}(\hat{x})\|^2\leq0,
\end{eqnarray}
then the desired $H_\infty$ filter for the system
(\ref{2017filter-6}) is given by (\ref{2017filter-7}).
\end{corollary}

\section{Suboptimal  $H_2/H_\infty$ filtering}

In this section, we further consider the suboptimal $H_2/H_\infty$
filtering design, that is, we design a filter that  not only
satisfies the robust $H_\infty$ performance, but also minimizes the
estimation error under the worst-case disturbance.

\begin{theorem}\label{2017filter-th6}
Consider system (\ref{2017filter-3}). For a given  disturbance
attenuation level $\gamma>0$, suppose there exist a positive
definite
 Lyapunov function sequence $\{V_k\}_{k\in{\mathcal N}}$, a
deterministic real-valued sequence $\{\gamma_k\geq0\}_{k\in{\mathcal
N}}$, and a nonnegative function $W: {\mathcal
R}^{n_{\eta}}\mapsto{\mathcal R}^+$ satisfying
(\ref{2017filter-rree}), (\ref{2017filter-L2}) and
\begin{equation}\label{2017filter-Ldd2dsd}
\Delta_{v} V_k(\eta)\leq \gamma_k-W(\eta),  \forall k\in {\mathcal
N}, \forall (\eta,v)\in {\mathcal R}^{n_\eta}\times {\mathcal
R}^{n_v},
\end{equation}
then the worst-case disturbance $\{v_k^*\}_{k\in {\mathcal N}}$ and
the corresponding augmented system state  $\{\eta_k^*\}_{k\in
{\mathcal N}}$ satisfy
\begin{equation}\label{gyggyb1}
H_k(\eta^*_k, v^*_k)-\gamma^2\|v^*_k\|^2=0, \ k\in {\mathcal N}.
\end{equation}
Moreover,
$$
\sum^{\infty}_{k=0}
E\|\tilde{z}^*_k\|^2=V_0(\eta_0)+\sum^{\infty}_{k=0} \gamma^2E\|{v}^*_k\|^2,
\ v=\{v_k^*\}_{k\in {\mathcal N}}.
$$
Moreover a suboptimal mixed $H_2/H_{\infty}$ filter can be
synthesized by solving the following constrained  optimization
problem:
\begin{equation*}
\min_{s.t.
\hat{f},\hat{g},\hat{m},(\ref{2017filter-rree}),(\ref{2017filter-Ldd2dsd})}
V_0(\eta_0).
\end{equation*}
\end{theorem}

\textbf {Proof}: Firstly, for any admissible external disturbance
$v=\{v_k\}_{k\in {\mathcal N}}\in l^2_{\infty}(\Omega, {\mathcal F},
{\{{\mathcal F}_k\}}_{k\in {\mathcal N}}, \mathcal{R}^{n_v})$
  and any initial state $\eta_0$,  by the smoothness of the conditional mathematical  expectation, we get
\begin{align}\label{ewdfc}
&EV_k(\eta_k)-EV_{k+1}(\eta_{k+1})\nonumber\\
=& E\{E[V_k(\eta_k)-V_{k+1}(\tilde{f}_k(\eta_k, w_k,v_k)]|{\mathcal
F}_k\}\nonumber\\
=&E\{[V_k(\eta)]-E[V_{k+1}(\tilde{f}_k(\eta, w_k,v)]\}_{\eta=\eta_k, v=v_k}\nonumber\\
=&E[-\Delta_v
V_k(\tilde{f}(\eta,w_k,v))+\gamma^2\|v\|^2-\|\tilde{z}_k\|^2]_{\eta=\eta_k,
v=v_k}-\gamma^2 E\|v_k\|^2+E\|\tilde{z}_k\|^2\nonumber\\
=&E[-H_k(\eta,v)+\gamma^2\|v\|^2]_{\eta=\eta_k, v=v_k}-\gamma^2 E\|v_k\|^2+E\|\tilde{z}_k\|^2.
\end{align}
Setting $v_k=v^*_k$ in (\ref{ewdfc}) and considering equation
(\ref{gyggyb1}), it follows that
$$
\gamma^2E\|v^*_k\|^2-E\|\tilde{z}^*_k\|^2=EV_{k+1}(\eta^*_{k+1})-EV_k(\eta^*_k),
$$
where $\tilde{z}^*_k$ is the estimation error corresponding to
$v_k^*$.  Taking the  summation from $k=0$ to $k=T$ on  both sides
of the above, we have
\begin{align}\label{eq vhvh}
\sum^T_{k=0}E[\gamma^2\|v^*_k\|^2-\|\tilde{z}^*_k\|^2]
=E[V_{T+1}(\eta^{*}_{T+1})]-V_0(\eta_0).
\end{align}
Obviously, for any admissible disturbance $\{v_k\}_{k\in {\mathcal
N}}\in l^2_{\infty}(\Omega, {\mathcal F}, {\{{\mathcal F}_k\}}_{k\in
{\mathcal N}}, \mathcal{R}^{n_v})$ satisfying
(\ref{2017filter-rree}), we have
\begin{align}\label{eqvhvhc}
\sum^T_{k=0}E[\gamma^2\|v_k\|^2-\|\tilde{z}_k\|^2]
\ge E[V_{T+1}(\eta_{T+1})]-V_0(\eta_0).
\end{align}
By (\ref{2017filter-Ldd2dsd}) and Lemma 3.2 of \cite{Lin2016}, it
yields that for $\forall\ v\in l^2_{\infty}(\Omega, {\mathcal F},
{\{{\mathcal F}_k\}}_{k\in \mathcal N}, \mathcal{R}^{n_v})$,
$$
\lim_{T\to\infty} E[V_{T+1}(\eta_{T+1})]=0,
$$
$$
\lim_{T\to\infty} E[V_{T+1}(\eta^*_{T+1})]=0.
$$
Letting $T\to\infty$ in (\ref{eq vhvh}) and (\ref{eqvhvhc}), we have
\begin{align*}
&\sum^{\infty}_{k=0}E[\gamma^2\|v_k\|^2-\|\tilde{z}_k\|^2]\\
\ge &\sum^\infty_{k=0}E[\gamma^2\|v^*_k\|^2-\|\tilde{z}^*_k\|^2]\\
=&-V_0(\eta_0),
\end{align*}
which shows that $\{v_k^*\}_{k\in {\mathcal N}}$ is the worst-case
disturbance, and
$$
\sum^{\infty}_{k=0}
E\|\tilde{z}^*_k\|^2=V_0(\eta_0)+\sum^{\infty}_{k=0} \gamma^2E\|{v}^*_k\|^2,
\ v=\{v_k^*\}_{k\in {\mathcal N}}.
$$
The theorem is proved.  $\square$

Based on Theorem \ref{2017filter-th6}, if we consider time-invariant
system (\ref{2017filter-new3})  and assume that $w_k$, $k\in
{\mathcal N}$,   have the same distribution, we can get the
following corollary that is  easily   verified.
\begin{corollary}\label{2017filter-cor1}
Consider system (\ref{2017filter-new3}). For a given  disturbance
attenuation level $\gamma>0$, suppose there exist a positive
 Lyapunov function $V(\eta)$, a
deterministic real-valued sequence $\{\gamma_k\geq0\}_{k\in{\mathcal
N}}$, and a nonnegative function $W: {\mathcal
R}^{n_{\eta}}\mapsto{\mathcal R}^+$, such that, for $\forall\ k\in
{\mathcal N}, \ \forall\ (\eta,v)\in {\mathcal R}^{n_\eta}\times
{\mathcal R}^{n_v}$,
\begin{equation}\label{2017fddilter-rreessw2}
H(\eta,v)-\gamma^2\|v\|^2\leq0,  \
\end{equation}
\begin{equation}\label{2017filter-Ldd2dsd2}
\Delta_{v} V(\eta)\leq \gamma_k-W(\eta),
\end{equation}
and
\begin{equation*}
\sum_{k=0}^{\infty}\gamma_k<\infty,
\end{equation*}
then the worst-case disturbance $\{v_k^*\}_{k\in {\mathcal N}}$ and
the corresponding augmented system state  $\{\eta_k^*\}_{k\in
{\mathcal N}}$ satisfy
\begin{equation*}
H(\eta^*_k, v^*_k)-\gamma^2\|v^*_k\|^2=0, \ k\in {\mathcal N}.
\end{equation*}
Moreover,
$$
\sum^{\infty}_{k=0}
E\|\tilde{z}^*_k\|^2=V(\eta_0)+\sum^{\infty}_{k=0} \gamma^2E\|{v}^*_k\|^2, \
v=\{v_k^*\}_{k\in {\mathcal N}}.
$$
Simultaneously, a suboptimal mixed $H_2/H_{\infty}$ filter can be
synthesized by solving the following constrained optimization
problem:
\begin{equation*}
\min_{s. t.
\hat{f},\hat{g},\hat{m},(\ref{2017fddilter-rreessw2}),(\ref{2017filter-Ldd2dsd2})}
V_0(\eta_0).
\end{equation*}
\end{corollary}

We find that for the general nonlinear stochastic system
(\ref{2017filter-1}), to design its mixed $H_2/H_{\infty}$ filter,
one needs to solve the constrained  optimization problem:
\begin{equation*}
\min_{s. t.
\hat{f},\hat{g},\hat{m},(\ref{2017filter-rree}),(\ref{2017filter-Ldd2dsd})}
V_0(\eta_0),
\end{equation*}
which is not an easy thing. However, for  linear discrete-time
stochastic systems, the above-mentioned problem can be converted
into solving a convex optimization problem. In particular, the
corresponding work for linear continuous time-invariant It\^o
systems has been done in \cite{Gershon}.

We consider the following system
\begin{eqnarray}\label{2017filter-linear1}
\begin{cases}
x_{k+1}=Ax_k+Bv_k+(Cx_k+Dv_k)w_k,\\
y_k=Kx_k+Lv_k,\\
z_k=Gx_k+Mv_k.
\end{cases}
\end{eqnarray}
Assume that $Ew_k=0$ and $Ew_k^2=1$ for all  $k\in {\mathcal N}$.
$\{w_k\}_{k\in {\mathcal N}}$ is the one-dimensional independent
random variable sequence.  We design the following filter for the
estimation of $z_k$:
\begin{eqnarray}\label{2017filter-linear2}
\begin{cases}
\hat{x}_{k+1}=Ax_k+\hat{K}(y_k-K\hat{x}_k),\\
\hat{z}_k=G\hat{x}_k.
\end{cases}
\end{eqnarray}
Denoting $\eta=\left[\begin{array}{ccc}x\\
x-\hat{x}\end{array}\right]$ and $\tilde{z}=z-\hat{z}$, we obtain
\begin{eqnarray}\label{2017filter-linear3}
\begin{cases}
\eta_{k+1}=\tilde{A}\eta_k+\tilde{B}v_k+(\tilde{C}\eta_k+\tilde{D}v_k)w_k,\\
\tilde{z}_k=\tilde{G}\eta_k+\tilde{M}v_k,
\end{cases}
\end{eqnarray}
where
\begin{eqnarray*}
\tilde{A}=\left[\begin{array}{cc}A&0\\
0&A-\hat{K}K\end{array}\right],\
\tilde{B}=\left[\begin{array}{cc}B\\ B-\hat{K}L\end{array}\right],\
\tilde{C}=\left[\begin{array}{cc}C&0\\ C&0\end{array}\right],\\
\tilde{D}=\left[\begin{array}{cc}D\\ D\end{array}\right],\
\tilde{G}=\left[\begin{array}{cc}0&G\end{array}\right],\
\tilde{M}=M.
\end{eqnarray*}
Setting $V(\eta)=\eta'P\eta$ and the disturbance attenuation level
$\gamma>0$, we have
\begin{align*}
&H(\eta,v)-\gamma^2\|v\|^2\\
=&E[\tilde{A}\eta+\tilde{B}v+(\tilde{C}\eta+\tilde{D}v)w_k]'P
[\tilde{A}\eta+\tilde{B}v\\
&+(\tilde{C}\eta+\tilde{D}v)w_k]-\eta'P\eta+(\tilde{G}\eta+\tilde{M}v)'\\
&(\tilde{G}\eta+\tilde{M}v)-\gamma^2v^2\\
=&\left[\begin{array}{ccc}\eta\\ v\end{array}\right]'\mathbb{P}
\left[\begin{array}{ccc}\eta\\ v\end{array}\right]\\
=&\eta'[\tilde{A}'P\tilde{A}+\tilde{C}'P\tilde{C}-P+\tilde{G}'\tilde{G}-
(\tilde{A}'P\tilde{B}+\tilde{C}'P\tilde{D}\\
&+\tilde{G}'\tilde{M})'(\tilde{B}'P\tilde{B}+\tilde{D}'P\tilde{D}+\tilde{M}'\tilde{M}-\gamma^2I)^{-1}(\tilde{A}'P\tilde{B}\\
&+\tilde{C}'P\tilde{D}+\tilde{G}'\tilde{M})]\eta
+(v-v^*)'(\tilde{B}'P\tilde{B}+\tilde{D}'P\tilde{D}\\
&+\tilde{M}'\tilde{M}-\gamma^2I)(v-v^*),
\end{align*}
where
\begin{align*}
\mathbb{P}=
\left[\begin{array}{ccc}\tilde{A}'P\tilde{A}+\tilde{C}'P\tilde{C}-P+\tilde{G}'\tilde{G}&\tilde{A}'P\tilde{B}+\tilde{C}'P\tilde{D}+\tilde{G}'\tilde{M}\\ *&
\tilde{B}'P\tilde{B}+\tilde{D}'P\tilde{D}+\tilde{M}'\tilde{M}-\gamma^2I\end{array}\right]
\end{align*}
and
\begin{eqnarray*}
v^*&=&-(\tilde{B}'P\tilde{B}+\tilde{D}'P\tilde{D}+\tilde{M}'\tilde{M}-\gamma^2 I)^{-1}\nonumber\\
&&\cdot(\tilde{A}'P\tilde{B}+\tilde{C}'P\tilde{D}+\tilde{G}'\tilde{M})\eta.
\end{eqnarray*}
Thus when the following  generalized algebraic Riccati inequality
(GARI)
\begin{eqnarray}\label{2017filter-ari1}
\begin{cases}
\tilde{A}'P\tilde{A}+\tilde{C}'P\tilde{C}-P+\tilde{G}'\tilde{G}-
(\tilde{A}'P\tilde{B}+\tilde{C}'P\tilde{D}\\+\tilde{G}'\tilde{M})'
(\tilde{B}'P\tilde{B}+\tilde{D}'P\tilde{D}+\tilde{M}'\tilde{M}-\gamma^2I)^{-1}\\
(\tilde{A}'P\tilde{B}+\tilde{C}'P\tilde{D}+\tilde{G}'\tilde{M})<0,\\
\tilde{B}'P\tilde{B}+\tilde{D}'P\tilde{D}+\tilde{M}'\tilde{M}-\gamma^2I<0
\end{cases}
\end{eqnarray}
admits  a positive definite matrix solution
$P>0$, $(\ref{2017filter-rree})$ holds. So system
(\ref{2017filter-linear3}) is  externally stable by
Lemma~\ref{2017filter-lemma2}.   From (\ref{2017filter-ari1}),
\begin{eqnarray*}
\Delta_{v=0}V(\eta)&=& \eta'(\tilde{A}'P\tilde{A}+\tilde{C}'P\tilde{C}-P)\eta\\
&<& -\eta'\Upsilon \eta\le 0,
\end{eqnarray*}
where,
\begin{eqnarray*}
\Upsilon&:=&\tilde{G}'\tilde{G}-
(\tilde{A}'P\tilde{B}+\tilde{C}'P\tilde{D}+\tilde{G}'\tilde{M})'\\
&&\cdot(\tilde{B}'P\tilde{B}+\tilde{D}'P\tilde{D}+\tilde{M}'\tilde{M}-\gamma^2I)^{-1}\\
&&\cdot(\tilde{A}'P\tilde{B}+\tilde{C}'P\tilde{D}+\tilde{G}'\tilde{M})\ge
0
\end{eqnarray*}
by
Lemmas~\ref{2017filter-lemmastablefffaa}-\ref{2017filter-lemmastablefff},
system (\ref{2017filter-linear3}) is internally stable.

Summarize the above discussions, we have

\begin{corollary}\label{2017filter-corlinear}
Consider system (\ref{2017filter-linear3}). For a given  disturbance
attenuation level $\gamma>0$, if  GARI (\ref{2017filter-ari1}) or
$\mathbb{P}<0$ has a positive definite matrix solution $P>0$, then
(\ref{2017filter-linear2}) is the $H_\infty$ filter of
(\ref{2017filter-linear1}). In this case, the worst-case disturbance
$\{v_k^*\}_{k\in {\mathcal N}}$ satisfies
\begin{eqnarray*}
v^*_k&=&-(\tilde{B}'P\tilde{B}+\tilde{D}'P\tilde{D}+\tilde{M}'\tilde{M}-\gamma^2 I)^{-1}\nonumber\\
&&\cdot(\tilde{A}'P\tilde{B}+\tilde{C}'P\tilde{D}+\tilde{G}'\tilde{M})\eta_k.
\end{eqnarray*}
Moreover, a suboptimal mixed $H_2/H_{\infty}$ filter can be
synthesized by solving the following constraint optimization
problem:
\begin{equation*}
\min_{s.t. \hat{K},(\ref{2017filter-ari1})} \ \text{trace}(P).
\end{equation*}
\end{corollary}

\begin{theorem}\label{2017filter-th7}
Consider system (\ref{2017filter-linear1}). For a given  disturbance
attenuation level $\gamma>0$, if there exist matrices $P_1>0$,
$P_2>0$ and $P_{K}$ solving the following LMI {\small
\begin{eqnarray}\label{2017filter-LMI1}
\left[\begin{array}{cccccccccccccccc}
-P_1&0&0&A'P_1&0&C'P_1&C'P_2&0\\
*&-P_2&0&0&A'P_2-K'P_K'&0&0&G\\
*&*&-\gamma^2I&B'P_1&B'P_2-L'P_K'&D'P_1&D'P_2&M'\\
*&*&*&-P_1&0&0&0&0\\
*&*&*&*&-P_2&0&0&0\\
*&*&*&*&*&-P_1&0&0\\
*&*&*&*&*&*&-P_2&0\\
*&*&*&*&*&*&*&-I
\end{array}\right]<0,\nonumber\\
\end{eqnarray}
} then the filter (\ref{2017filter-linear2}) is the desired
$H_{\infty}$ filter and the  filter parameter is given by
$$
\hat{K}=P_2^{-1}P_K.
$$
Moreover, a suboptimal mixed $H_2/H_{\infty}$ filter can be
synthesized by solving the following constraint optimization
problem:
\begin{equation*}
\min_{s.t. \hat{K},(\ref{2017filter-LMI1})}\ \text{trace}(P_1+P_2).
\end{equation*}
\end{theorem}

\textbf {Proof}: By Schur's complement, $\mathbb{P}<0$ is equivalent
to that
\begin{equation*}
\left[\begin{array}{ccccccccccccccccc}
-P&0&\tilde{A}'P&\tilde{C}'P&\tilde{G}'\\
*&-\gamma^2I&\tilde{B}'P&\tilde{D}'P&\tilde{M}'\\
*&*&-P&0&0\\
*&*&*&-P&0\\
*&*&*&*&-I
\end{array}\right]<0.
\end{equation*}
If we take $P=diag(P_1,P_2)$, considering  system
(\ref{2017filter-linear3}), then this theorem is proved
 by Corollary \ref{2017filter-corlinear}.

\section{Illustrative examples}
In this section, we give two examples to illustrate the
effectiveness of our obtained results.

\begin{figure}
  \centering
  \includegraphics[width=3.5in]{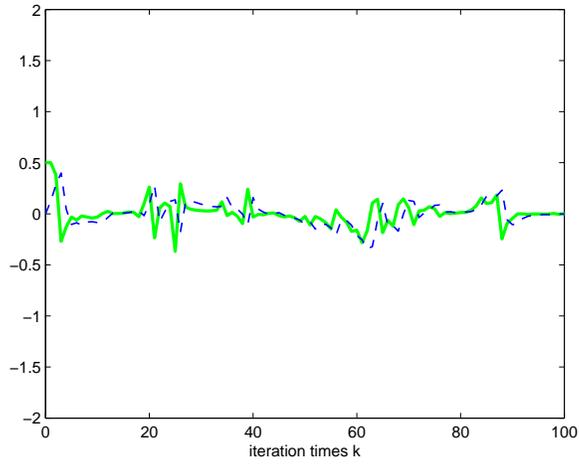}
  \caption{$z_k$ and $\hat{z}_k$.}
  \label{fig1}
\end{figure}

\begin{figure}
  \centering
  \includegraphics[width=3.5in]{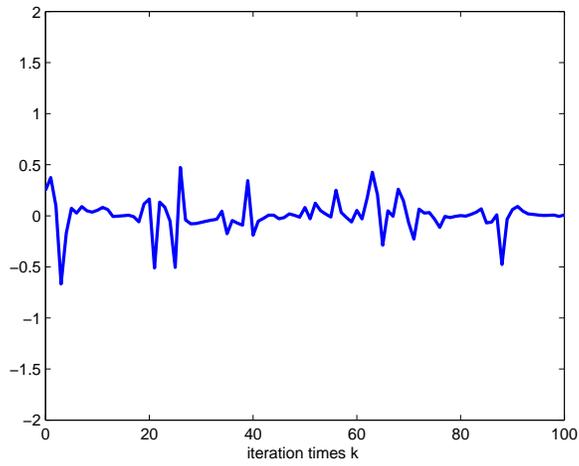}
  \caption{The estimation error $\tilde{z}_k$.}
  \label{fig2}
\end{figure}

\begin{figure}
  \centering
  \includegraphics[width=3.5in]{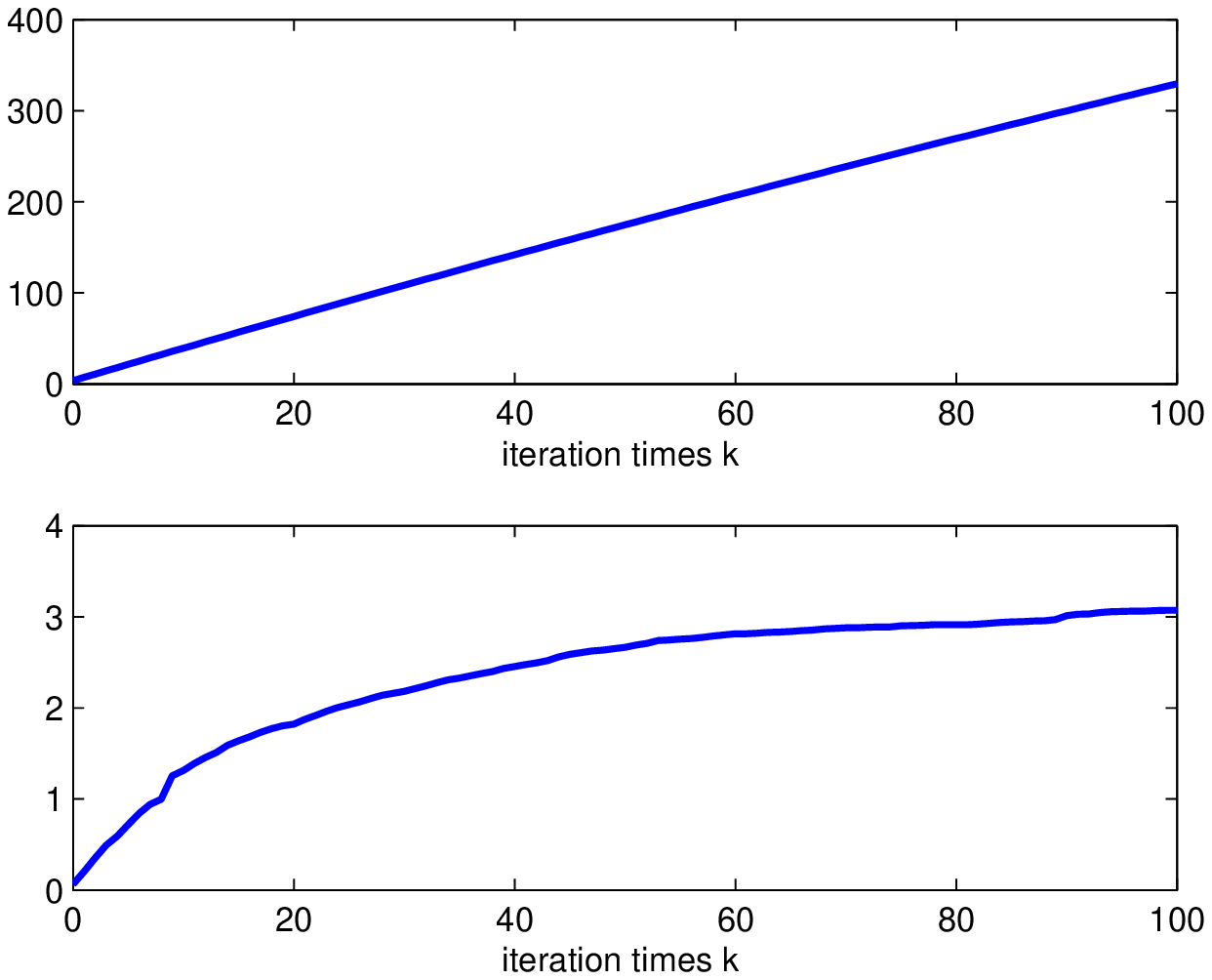}
  \caption{$\sum_{i=0}^{k}E\gamma^2v_i^2$ and $\sum_{i=0}^{k}E\tilde{z}_i^2$.}
  \label{fig3}
\end{figure}

\textbf {Example 5.1}: Let
$x_k=\left[\begin{array}{ccc}x_k(1)\\x_k(2)\end{array}\right]$,
$y_k=\left[\begin{array}{ccc}y_k(1)\\y_k(2)\end{array}\right]$. We
consider the following nonlinear discrete-time stochastic system:
\begin{equation}\label{2017filter-example1}
\begin{cases}
 \  x_{k+1}(1)=\displaystyle\frac{0.6x_{k}(1)^3}{1+x_{k}(1)^2+x_{k}(2)^2}+0.1v_kx_{k}(1)\\
 \ \ \ \ \ \ \ \ \ \ \ \ \ \ \ \ +0.5x_{k}(2)sin(v_k)w_k, \\
 x_{k+1}(2)=0.65x_{k}(2)+0.1x_{k}(2)x_{k}(1)\\
 \ \ \ \ \ \ \ \ \ \ \ \ \ \ \ \ +0.5v_ksin(x_{k}(2))w_k, \\
 \  y_k(1)=0.5x_k(1)+v_ksin(x_k(1)),\\
 y_k(2)=0.5x_k(2)+v_ksin(x_k(2)),\\
    z_k=0.1x_k(1)+0.1x_k(2).
   \end{cases}
\end{equation}
In the sequel, we design an $H_{\infty}$ filter for system
(\ref{2017filter-example1}). The environmental noise $\{w_k\}_{k\in
{\mathcal N}}$ is assumed as a one-dimensional independent  white
noise process and the external disturbance $v_k=2\times(0.9999)^k$,
$k\in {\mathcal N}$, so $\sum_{k=0}^{\infty}v_k^2<\infty$. We choose
the Lyapunov function as
$$
V(\eta_k)=\|x_k\|^2+\|\hat{x}_k\|^2= \left[\begin{array}{ccc}
x_k\\
{\hat x}_k
\end{array}
\right]' \left[\begin{array}{ccc}
I & 0\\
0 & I
\end{array}
\right]\left[\begin{array}{ccc}
x_k\\
{\hat x}_k
\end{array}
\right].
$$
By Corollary \ref{2017filter-cor1}, it is easy to test  that
(\ref{eq44}) and (\ref{eq45}) hold with $Q_1=I$ and $Q_2=I$ and the
appropriate $H_{\infty}$ filter for (\ref{2017filter-example1}) is
designed as
\begin{equation}\label{2017filter-example2}
\begin{cases}
 \hat{x}_{k+1}(1)=0.5\hat{x}_k(1)+0.5y_k(1),\\
  \hat{x}_{k+1}(2)=0.5\hat{x}_k(2)+0.5y_k(2),\\
 \hat{z}_k=0.1\hat{x}_k(1)+0.1\hat{x}_k(2).
   \end{cases}
\end{equation}
Figures \ref{fig1} and \ref{fig2} show a sample trajectory in an
experiment. Figure \ref{fig1} shows trajectories  of the $z_k$ and
$\hat{z}_k$, and the error $\tilde{z}_k=z_k-\hat{z}_k$ is depicted
in Figure \ref{fig2}. We use Matlab to simulate system
(\ref{2017filter-example1}) and system (\ref{2017filter-example2})
for 1000 times to obtain the approximate value of
$\sum_{i=0}^{k}E\tilde{z}_i^2$. In Figure \ref{fig3}, we can see
that $\sum_{i=0}^{k}E\tilde{z}_i^2$ is always less than
$\sum_{i=0}^{k}E\gamma^2v_i^2$, which is in accordance with our
theoretical analysis.

\begin{figure}
  \centering
  \includegraphics[width=3.5in]{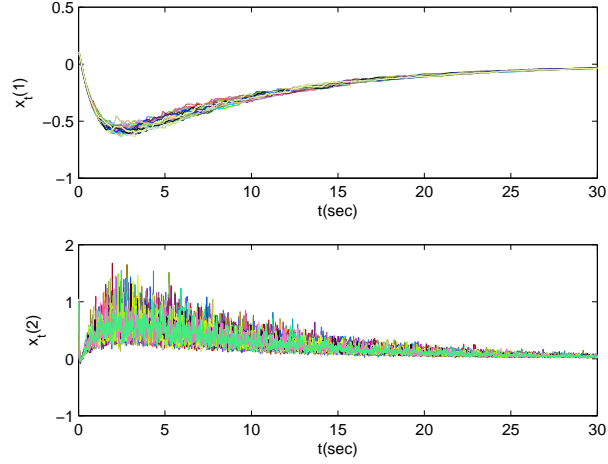}
  \caption{The state trajectory $x_k$ of system (\ref{2017filter-system9}).}
  \label{x}
\end{figure}

\begin{figure}
  \centering
  \includegraphics[width=3.5in]{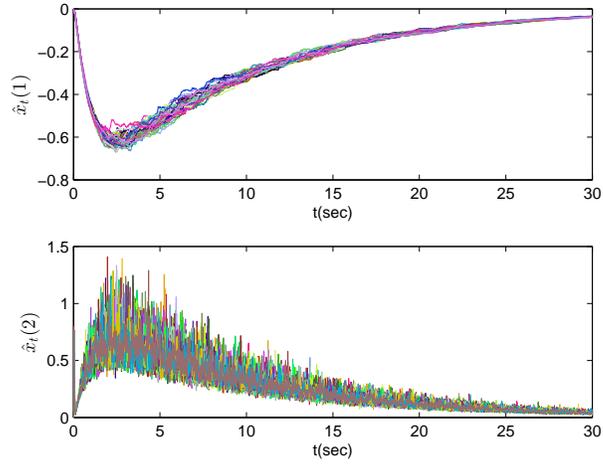}
  \caption{The state trajectory $\hat{x}_k$ of system (\ref{2017filter-system10}).}
  \label{hatx}
\end{figure}

\begin{figure}
  \centering
  \includegraphics[width=3.5in]{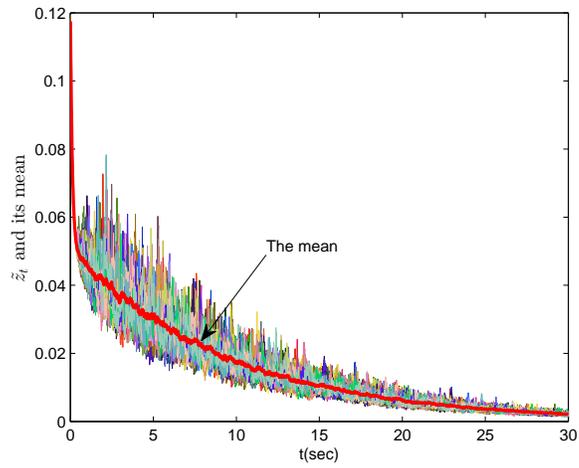}
  \caption{$\tilde{z}_k(1)=x_k(1)-\hat{x}_k(1)$.}
  \label{tildez}
\end{figure}

\begin{figure}
  \centering
  \includegraphics[width=3.5in]{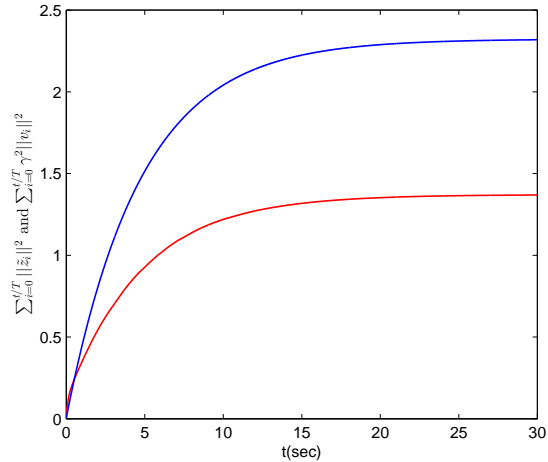}
  \caption{$\sum^{k}_{i=0}E\|\tilde{z}_{i}\|^2$ and $\sum^{k}_{i=0}E\gamma^2\|{v}_{i}\|^2$.}
  \label{z and v}
\end{figure}

\textbf {Example 5.2}: In order to verify the validity of Theorem
\ref{2017filter-th7}, we use the vehicle model in $[1]$. A vehicle's roll dynamic is governed by the following
differential equation:
\begin{align}\label{eqvhvhrrr}
&I_{xx}\ddot{\eta}+(C_R+D_nw_k)\dot{\eta}+K_R\eta\nonumber\\
=&m_sa_yh_{cr}+m_sh_{cr}g\sin(\eta),
\end{align}
where  $\eta$ is the vehicle roll angle,  $D_n$ means the noise
intensity, $I_{xx}$ is the sprung mass moment of the inertia with
respect to the roll axis, $m_s$ is the sprung mass, $h_{cr}$ is the
sprung mass height about the roll axis, $C_R$ is the total torsional
damping, $K_R$ is the stiffness coefficient, $a_y$ is the lateral
acceleration at the vehicle center of gravity (COG) and $g$ is the
acceleration due to gravity.  $\{w_k\}_{k\in {\mathcal N}}$
represents the  system internal  noise driven by one-dimensional
independent white noise processes  with $E[w_k]=0$,
$E[w_kw_j]=\delta_{kj}$,  where   $\delta_{kj}$ is   a Kronecker
function defined by $\delta_{kj}=0$ for $k\ne j$ while
$\delta_{kj}=1$ for $k=j$. Then, by  setting the length of the
sampling interval $T=0.01$, the continuous-time system
(\ref{eqvhvhrrr}) can be discretized into the following system:
\begin{align}\label{2017filter-system9}
\begin{cases}
x_{k+1}=Ax_k+Cx_kw_k+Bv_k,\\
y_{k}=Hx_k+Lv_k,\\
z_k=Mx_k,
\end{cases}
\end{align}
where
$x_k=\left[\begin{array}{ccc}x_k(1)\\x_k(2)\end{array}\right]=\left[\begin{array}{ccc}\eta_k\\
\Delta \eta_k\end{array}\right]$,  $v_k$ is the disturbance with
$v_k=0.01\times0.9^{t}$, $y_{k}$ is the measurement signal, $z_k$ is
the regulation output, and
\begin{align*}
&A=\left[\begin{array}{cccc} 1&T\\
\frac{(m_sh_{cr}-K_R)T}{I_{xx}}&1-\frac{C_RT}{I_{xx}}
\end{array}\right],\
C=\left[\begin{array}{cccc}0&0\\0&-\frac{D_nT}{I_{xx}}
\end{array}\right],\\
&H=\left[\begin{array}{cccc} 0.4048&0.6405\\
    1.1213&1.4616\end{array}\right],\
B=\left[\begin{array}{cccc}  -0.7916\\
    0.3652\end{array}\right],\\
&L=\left[\begin{array}{cccc} 0.8248\\
   -1.3774\end{array}\right],\
M=\left[\begin{array}{cccc}1&0\end{array}\right].
\end{align*}
All parameters of the  Mercedes-Benz commercial vehicle used in
$[1]$ are presented in Table 1.
\begin{table}[!htb]\label{2018filter-table1}\caption{}
\centering \scalebox{1}{
\begin{tabular}{c|c|c}
\hline
\multicolumn{3}{c}{Parameters of the  Mercedes-Benz commercial vehicle}\\
 \hline
  Symbol&Value&Unit\\
  \hline
  % \> for next tab, \\ for new line...
  $C_R$&53071& N ms/rad \\
  $m_s$&1700 &kg \\
  $h_{cr}$&0.25& m \\
  $I_{xx}$& 1700 &kg $m^2$ \\
  $K_R$&55314& N ms/rad\\
  $D_n$&20& N ms/rad\\
  \hline
  \end{tabular}}
  \end{table}
Then, in order to estimate $\eta_k$, we need to determine the
parameter $\hat{K}$. So, by Theorem \ref{2017filter-th7},  we can
find a set of feasible solutions to (\ref{2017filter-LMI1}) as
follows:
\begin{align*}
P_1=\left[\begin{array}{cccc}0.0114&0.0002\\
    0.0002&0.0002\end{array}\right],\ \\
P_2=\left[\begin{array}{cccc} 7.5939&0.1379\\
    0.1379&0.0029\end{array}\right],\\
P_K=\left[\begin{array}{cccc}-6.3529&2.8009\\
   -0.1137&0.0503\end{array}\right]
\end{align*}
and
\begin{align*}
\hat{K}=\left[\begin{array}{cccc} -0.9194&0.3965\\
    4.5617&-1.5241\end{array}\right].
\end{align*}
Thus we can design a proper filter as
\begin{align}\label{2017filter-system10}
\begin{cases}
\hat{x}_{k+1}=\left[\begin{array}{cccc} 1&T\\
\frac{(m_sh_{cr}-K_R)T}{I_{xx}}&1-\frac{C_RT}{I_{xx}}
\end{array}\right]\hat{x}_{k}\\
\ \ \ \ \ \ \ \ \ \ \ \ +\left[\begin{array}{cccc} -0.9194&0.3965\\
    4.5617&-1.5241\end{array}\right](y_k-H\hat{x}_k),\\
   \hat{z}_k=\hat{x}_k(1)
\end{cases}
\end{align}
under $\hat{x}_0=\left[\begin{array}{cccc} \hat{x}_0(1)\\
   \hat{x}_0(2)\end{array}\right]=\left[\begin{array}{cccc} 0\\
   0\end{array}\right]$.
Using Matlab to simulate systems
(\ref{2017filter-system9})-(\ref{2017filter-system10}) for 100
times under $x_0=\left[\begin{array}{cccc} x_0(1)\\
   x_0(2)\end{array}\right]=\left[\begin{array}{cccc} 0.1\\
   1\end{array}\right]$, we can obtain  Figures \ref{x}-\ref{hatx}. From Figures \ref{x} and \ref{hatx}, we can see that the augmented
system is stable. Figure \ref{tildez} displays
$\tilde{z}_k(1)=x_k(1)-\hat{x}_k(1)$, which  converges to zero
quickly. So, filter (\ref{2017filter-system10}) can track the the
adjustment output of (\ref{2017filter-system9}). Based on the data
of the 100 experiments, we obtain the approximate value of
$\sum^{k}_{i=0}E\|\tilde{z}_{i}\|^2$. In Figure \ref{z and v}, the
red curve stands for $\sum^{k}_{i=0}E\|\tilde{z}_{i}\|^2$ and blue
curve represents $\sum^{k}_{i=0}E\gamma^2\|v_{i}\|^2$. Figure \ref{z
and v} shows that $\sum^{k}_{i=0}E\|\tilde{z}_{i}\|^2
\leq\sum^{k}_{i=0}\gamma^2 E\|v_{i}\|^2$,
 which is in accordance with our
theoretical analysis.


\section{Conclusions\label{sec:Concl}}

\hspace{0.13in} This paper has  studied  the robust $H_\infty$
filtering  of general nonlinear  discrete  stochastic systems. A
SBRL has been obtained based on the property of a conditional
mathematical expectation (Lemma 2.2). By means of the discrete-time
stochastic LaSalle's invariance principle,  it is shown that the
nonlinear stochastic $H_\infty$ filtering can be constructed by
solving an  HJI.  In the case of the worst-case disturbance
$\{v^*_k\}_{k\in {\mathcal N}}$, a suboptimal $H_2/H_\infty$
filtering has also been  studied. Two examples including a practical
example are presented to illustrate the validity of our main
results.

\end{document}